\documentclass[12pt]{article}
\usepackage{bbm}
\usepackage{mathrsfs}
\usepackage{amsfonts}
\usepackage{amsmath}
\usepackage{amssymb,amsmath}
\def\Z{\mathbb{Z}}
\def\C{\mathbb{C}}
\def\cl{\centerline}
\def\vs{\vspace*}

\def\QED{\hfill$\Box$}

\def\D{\Delta}
\def\A{\mathcal{A}}
\def\a{\alpha}
\def\d{\cdot}
\def\c{\circ}
\def\o{\otimes}
\def\p{\partial}
\def\l{\lambda}
\def\L{\mathcal{L}(\mathcal{A})}
\def\LL{\mathcal{L}}
\def\dis{\displaystyle}

\numberwithin{equation}{section}
\newtheorem{theo}{Theorem}[section]
\newtheorem{defi}[theo]{Definition}
\newtheorem{exam}[theo]{Example}
\newtheorem{coro}[theo]{Corollary}

\newtheorem{prop}[theo]{Proposition}

\newtheorem{rema}[theo]{Remark}

\def\dis{\displaystyle}

\begin{document}

\cl{{\large\bf  Hom Gel'fand-Dorfman bialgebras and Hom-Lie
conformal algebras} \footnote{Corresponding author:
lmyuan@hit.edu.cn}} \vskip8pt

\cl{ Lamei Yuan} \vskip1mm
 \cl{\small{Science Research Center, Academy of Fundamental and Interdisciplinary
 Sciences,}}
\cl{\small{Harbin Institute of Technology, Harbin 150080, China}}
\cl{\small E-mail: lmyuan@hit.edu.cn
 }\vs{6pt}

{\small\parskip .005 truein \baselineskip 3pt \lineskip 3pt

\noindent{\bf Abstract.}  The aim of this paper is to introduce the notions of
Hom Gel'fand-Dorfman bialgebra and Hom-Lie conformal algebra. In this paper, we give four constructions of Hom Gel'fand-Dorfman bialgebras. Also, we provide a general construction of Hom-Lie conformal algebras from Hom-Lie algebras. Finally, we prove that a Hom Gel'fand-Dorfman bialgebra is equivalent to a Hom-Lie conformal algebra of degree $2$.
 \vs{5pt}


\parskip .001 truein\baselineskip 6pt \lineskip 6pt

\vs{10pt}

\cl{\bf\S1. \
Introduction}\setcounter{section}{1}\setcounter{equation}{0}
\vs{6pt}

The notion of Novikov algebra was originally introduced in
connection with the Poisson bracket of hydrodynamic type
\cite{BN,DN1,DN2} and Hamiltonian operators in formal
variational calculus \cite{GD,Xu1}. The systemetic study of Novikov
algebras was started by Zel'manov \cite{Z} and Filipov \cite{F},
whereas the term ``Novikov algebra" was first used by Osborn
\cite{O1}. Novikov algebras constitute a special class of left-symmetric algebras (or
under other names like pre-Lie algebras, Vinberg algebras and quasi-associative algebras), arising from the study of affine manifolds, affine structures and convex homogeneous cones
\cite{A,K,V}. Left-symmetric algebras are closely related with many
important fields in mathematics and mathematical physics, such as
infinite-dimensional Lie algebras \cite{BM1}, classical and
quantum Yang-Baxter equation \cite{ES,GS}, and quantum field theory
\cite{CK}. The superanalogue of Novikov algebras was introduced and studied by Xu \cite{Xu3, Xu4}.

Hom-Lie algebras were initially introduced by Hartwig, Larsson and Silvestrov \cite{HLS} motivated by constructing examples of deformed Lie algebras coming from twisted discretizations of vector fields. This kind of algebraic structure contains Lie algebras as a particular case. Hom-Lie algebras were intensively discussed in \cite{Yau1,Yau2,Yau3,Yau4}. Hom-associative algebras were introduced in \cite{MS}, where it was shown that the commutator
bracket of a Hom-associative algebra gives rise to a Hom-Lie algebra. The graded case was mentioned in \cite{LS1,LS2}, while a detailed study was given in \cite{AM,Yuan}.

Gel¡¯fand and Dorfman \cite{GD} introduced the notion of
Gel'fand-Dorfman bialgebra in the study of certain Hamiltonian
pairs, which play important roles in complete integrability of
nonlinear evolution partially differential equations. The
super Gel'fand-Dorfman bialgebra was introduced and studied in
\cite{X1}, where it was proved that a super Gel'fand-Dorfman bialgebra is equivalent to a quadratic conformal superalgebra. An
analogous result for Gel'fand-Dorfman bialgebras was essentially
stated in \cite{GD}.

The notion of Lie conformal (super)algebra, introduced by
Kac \cite{ka}, encodes an axiomatic description of the operator
product expansions of chiral fields in conformal field theory. Closely related to vertex
algebras \cite{Li1,Li2}, Lie conformal algebras have important applications in other areas of algebras and
integrable systems \cite{Z1,Z2}. In particular, they provide
powerful tools in the study of infinite-dimensional Lie
(super)algebras satisfying the locality property \cite{ka,SY}. The main examples are those based on the punctured complex plane, such as the Virasoro algebra and the loop Lie algebras \cite{DK}.

In this paper, we introduce and study Hom Gel'fand-Dorfman bialgebras and Hom-Lie conformal algebras, which are natural generalizations of Gel'fand-Dorfman bialgebras and Lie conformal algebras, respectively.

 The paper is organized as follows: In Sec. 2, we review the notions of Novikov algebra, Hom-Lie algbra, Hom-Novikov algebra, Gel'fand-Dorfman algebra and Lie conformal algebra. In Sec. 3, we introduce the notion of Hom Gel'fand-Dorfman bialgebra and give four basic constructions of Hom Gel'fand-Dorfman bialgebras. In Sec. 4, we introduce the notion of Hom-Lie conformal algebra
 and provide a general construction of Hom-Lie conformal algebras. In Sec. 5,
we prove that a Hom Gel'fand-Dorfman bialgebra is equivalent to a Hom-Lie conformal algebra of degree 2.
\vs{12pt}

\cl{\bf\S2. \
Preliminaries}\setcounter{section}{2}\setcounter{equation}{0}
\vs{6pt}
In this section, we recall the notions of Novikov algebra, Hom-Novikov algebra, Hom-Lie algbra, Gel'fand-Dorfman algebra and Lie conformal algebra. Throughout this paper, all vector spaces and tensor products are over the complex field $\mathbb{C}$. In addition to the standard notation $\Z$, we use $\Z^+$ to denote the set of nonnegative
integers.
\begin{defi}{\rm
A Novikov algebra is a vector space $\A$ equipped with
an operation $\circ$ such that for $\ x, y,z \in\A$:
\begin{eqnarray}\label{Nov11}
(x\circ y)\circ z =(x\circ z)\circ y,\ \ (x,y,z)= (y,x,z),
\end{eqnarray}
where the associator $(x,y,z)=(x\circ y)\circ z-x\circ (y\circ z).$ }
\end{defi}

The first identity in
\eqref{Nov11} says that the right multiplication operators on $\A$
commute with one another and the second identity states that the associator is
symmetric in the first two variables.

Yau introduced a twisted generalization of Novikov
algebras, named Hom-Novikov algebras \cite{Yau}.

\begin{defi}{\rm  A Hom-Novikov algebra is a vector space $\mathcal {A}$ equipped with a bilinear
operation $\c$ and a linear endomorphism $\alpha$, such that the
following two relations hold for $x,y,z\in\mathcal{A}$:
\begin{eqnarray}
(x\c y)\c\alpha(z)-\alpha(x)\c(y\c z)&=&(y\c x)\c \alpha(z)
-\alpha(y)\c (x\c z),\label{Hom-Nov1+}\\
(x\c y)\c \alpha(z)&=&(x\c z)\c \alpha(y).\label{Hom-Nov2+}
\end{eqnarray}
}\end{defi}

Obviously, Novikov algebras are examples of Hom-Novikov algebras where $\a={\rm id}$. The study of Hom-type algebras started with the notion of Hom-Lie algebra, which was originally
introduced by Hartwig, Larsson and Silvestrov \cite{HLS} motivated by constructing deformations of the Witt and Virasoro algebras based on $\sigma$-derivations.

\begin{defi} {\rm A Hom-Lie algebra is a vector space $L$ with a
bilinear  map  $[ \cdot,\cdot]:L\times L\longrightarrow L$ and a
linear map $\alpha:L\longrightarrow L$, such that the following relations hold for $x,y,z\in L$:
\begin{eqnarray}
&&[x,y]=-[y,x],\hspace{0.5cm}\mbox{(skew-symmetry)}\label{ss}\\
&&[[x,y],\alpha(z)]+[[y,z],\alpha(x)]+[[z,x],\alpha(y)]=0.\hspace{0.5cm}\mbox{(Hom-Jacobi
\  identity)}\label{JJ}\ \ \ \ \ \ \
\end{eqnarray}}
\end{defi}

The notion of Gel'fand-Dorfman bialgebra was initially introduced in the study of Hamilton operators \cite{GD}.
\begin{defi}{\rm A  Gel'fand-Dorfman bialgebra is a vector space $\A$, equipped with two operations
$[\cdot, \cdot]$ and $\circ$ such that $(\A, [\cdot, \cdot])$ forms a Lie
algebra, $(\A, \circ)$ forms a Novikov algebra and the compatibility
condition holds for $x, y, z\in\A$:
\begin{eqnarray}\label{GB}
[x\circ y, z]-[x\circ z, y]+[x,y]\circ z-[x,z]\circ y-x\circ[y,z]=0.
\end{eqnarray}}
\end{defi}

Such a bialgebraic structure corresponds
to the following Poisson bracket of dynamic type
\begin{eqnarray}\label{diss}
[u(x),v(y)]= [u, v](x)\delta(x,y)+\partial_x(u\c v)(x)\delta(x,
y)+(u\c v+v\c u)(x)\partial_x\delta(x,y).
\end{eqnarray}
It was pointed out in \cite{GD} that if we define
\begin{eqnarray}\label{comm}
[x,y]^-=x\c y-y\c x, \ \mbox{for}\ x,y\in\A,
\end{eqnarray} for a Novikov algebra $(\A, \circ)$, then $(\A, [\cdot,
\cdot]^-,\circ)$ forms a Gel'fand-Dorfman bialgebra. 

 The follow definition is due to Kac \cite{ka}.
\begin{defi}\label{definition}{\rm
 A Lie conformal algebra is a $\C[\partial ]$-module $\mathcal {R}$ with a
$\lambda$-bracket $[\cdot_\lambda\cdot]$ which defines a
$\C$-bilinear map from $\mathcal {R}\otimes \mathcal {R}$ to
$\mathcal {R}[\lambda]$, where $\mathcal
{R}[\lambda]=\C[\lambda]\otimes \mathcal {R}$ is the space of
polynomials of $\lambda$ with coefficients in $\mathcal {R}$,
satisfying the following axioms for all $a, b, c\in \mathcal {R}$:
\begin{eqnarray}
[\partial a_\lambda b]&=&-\lambda[a_\lambda b],\ \ \ \ \ \mbox{(conformal\  sesquilinearity)}\label{L1}\\
{[a_\lambda b]} &=& -[b_{-\lambda-\partial}a],\ \ \ \ \ \ \mbox{(skew-symmetry)}\label{L2}\\
{[a_\lambda[b_\mu c]]}&=&[[a_\lambda b]_{\lambda+\mu
}c]+[b_\mu[a_\lambda c]].\ \ \mbox{(Jacobi \ identity)}\label{L3}
\end{eqnarray}
}\end{defi}

Let $m$ be a positive integer. A
Lie conformal algebra $\mathcal {R}$ is called a {\it Lie conformal algebra of degree $m$} if for any $a, b \in \mathcal
{R}$,
\begin{eqnarray}
[a_\l b]=\mbox{$\sum\limits_{j\geq0,i+j< m}$}\p^i w_{i,j}\l^{j}, \ \mbox{for\ some} \ \ w_{i,j}\in\mathcal {R}.
\end{eqnarray}

It was stated in \cite{GD} that a Gel'fand-Dorfman bialgebra is equivalent to a Lie conformal algebra of degree $2$.  This statement was proved by Xu \cite{X1} in the super case.
\vs{10pt}

\cl{\bf\S3. Hom Gel'fand-Dorfman
bialgebras}\setcounter{section}{3}\setcounter{equation}{0}\setcounter{theo}{0}
\vs{8pt}

In this section, we introduce the notion of Hom Gel'fand-Dorfman bialgebra. Also, we present four constructions of Hom Gel'fand-Dorfman bialgebras from Hom-Novikov algebras, Gel'fand-Dorfman bialgebras, commutative Hom-associative algebras and Hom-Poisson algebras along
with some suitable algebra endomorphisms or derivations.

\begin{defi}{\rm  A Hom
Gel'fand-Dorfman
 bialgebra is a vector space $\A$ equipped with a linear endomorphism $\a$ and two operations
 $[\cdot, \cdot]$ and $\circ$, such that $(\A,
[\cdot, \cdot], \a)$ is a Hom-Lie algebra, $(\A, \circ, \a)$
is a Hom-Novikov algebra and the following compatibility condition holds for $x, y, z\in\A$:
\begin{eqnarray}\label{HGB}
[x\circ y, \a(z)]-[x\circ z, \a(y)]+[x,y]\circ \a(z)-[x,z]\circ
\a(y)-\a(x)\circ[y,z]=0.
\end{eqnarray}}
\end{defi}

Clearly, we recover the Gel'fand-Dorfman bialgebras when $\a={\rm id}$. In the sequel, we shall simply call a  Hom Gel'fand-Dorfman bialgebra a
 Hom-GD bialgebra.

Our first construction of Hom-GD bialgebras is from Hom-Novikov algebras, analogous to
the fundamental construction of Gel'fand-Dorfman bialgebras from Novikov algebras via the commutator bracket \cite{GD}.

\begin{theo} \label{th1} Let $(\A,\circ,\a)$ be a Hom-Novikov algebra. Define the commutator
\begin{eqnarray}\label{commu**}
[x,y]^-=x\c y-y\c x, \  \ \mbox{for} \ x, y\in\A.
\end{eqnarray}
Then $(\mathcal {A},[\cdot,\cdot]^-, \circ, \alpha)$ forms a Hom-GD algebra.
\end{theo}
\noindent{\it Proof.~} Clearly, $[\cdot,\cdot]^-$ is skew-symmetric. For any $x,y,z\in\A$, we have \begin{eqnarray}\label{j1}
&&[x,y]^-\c \alpha(z)+[y,z]^-\c \alpha(x)+[z,x]^-\c \alpha(y)\nonumber\\
&=&\big((x\c y)\c \alpha(z)-(x\c z)\c \alpha(y)\big)+\big((y\c z)\c \alpha(x)-(y\c x)\c \alpha(z)\big)\nonumber\\&&+\big((z\c x)\c \alpha(y)-(z\c y)\c \alpha(x)\big)\nonumber\\
&=&0,\end{eqnarray}
by \eqref{Hom-Nov2+} and \eqref{commu**}. Similarly, we get
\begin{eqnarray}\label{j2}
\alpha(x)\c [y,z]+\alpha(y)\c [z,x]+\alpha(z)\c [x,y]=0,
\end{eqnarray}
which together with \eqref{j1} implies \eqref{JJ}. Thus $(\A, [\cdot,\cdot]^-, \alpha)$ is a Hom-Lie algebra. By \eqref{Hom-Nov1+} and \eqref{commu**}, we obtain
\begin{eqnarray*}
&&{[x \circ y, \a(z)]^{-} -[x\circ z, \a(y)]^{-} +[x,y]^{-}\circ
\a(z)-[x,z]^{-} \circ \a(y)- \a(x)\circ [y,z]^{-}}\\
&=&(x \circ y)\circ \a(z)- \a(z)\circ(x \circ y)-(x\circ
z)\circ\a(y)+\a(y)\circ(x\circ z)+(x\circ y)\circ\a(z)\\&& -(y\circ
x)\circ\a(z)-(x\circ z)\circ\a(y)+(z\circ
x)\circ\a(y)-\a(x)\circ(y\circ z)+\a(x)\circ(z\circ y)\\
&=&\big((x\circ y)\circ\a(z) +\a(y)\circ(x\circ z)-(y\circ
x)\circ\a(z)- \a(x)\circ(y \circ z)\big)\\&& +\big((z\circ
x)\circ\a(y)+\a(x)\circ(z\circ y)- \a(z)\circ(x \circ y)-(x\circ
z)\circ\a(y)\big)\\&=&0.
\end{eqnarray*}  
This proves \eqref{HGB} and the theorem.
\QED\vskip8pt

 The following result, due to Yau \cite{Yau}, gives a construction of Hom-Novikov algebras
 from Novikov algebras with an algebra endomorphism.
\begin{prop} \label{prop2} Let $(\mathcal
{A}, \circ)$  be a Novikov algebra with an algebra
endomorphism $\a$. Then $(\mathcal {A}, \circ_{\a}, \a)$ forms a
Hom-Novikov algebra, where \begin{eqnarray}\label{****}
x \circ_{\a} y =\alpha(x)\circ
\alpha(y),\ \mbox{ for}\  x,y\in\A.
\end{eqnarray}
\end{prop}

 As a consequence of Theorem \ref{th1} and Proposition
\ref{prop2}, we obtain
\begin{coro}\label{co} Let $(\mathcal
{A}, \circ)$  be a Novikov algebra and $\a$ an algebra endomorphism of $(\mathcal
{A}, \circ)$. Then $(\mathcal {A}, [\cdot,\cdot]^-_{\a},\circ_{\a}, \a)$ is a Hom-GD bialgebra,  with $\circ_{\a}$ defined by \eqref{****} and
\begin{eqnarray}
[x,y]^-_{\a}=\a(x\circ y)-\a(y\circ x),\ \mbox{for}\ x,y\in\A.
\end{eqnarray}
\end{coro}

The following theorem, which is an analogue of Proposition \ref{prop2}, shows our second construction of Hom-GD bialgebras from a Gel'fand-Dorfman algebra with an algebra endomorphism.
\begin{theo}\label{p-2}
Let $(\mathcal {A}, [\cdot,\cdot],\circ)$ be a Gel'fand-Dorfman
algebra with an algebra endomorphism $\a$. Then $(\mathcal {A},
[\cdot,\cdot]_\a,\circ_\a, \a)$ forms a Hom-GD bialgebra, where $\circ_{\a}$ is defined by \eqref{****} and 
\begin{eqnarray}[x,y]_\a=[\a(x),\a(y)], \ \mbox{for}\ x,y\in\A.\end{eqnarray}
\end{theo}
\noindent{\it Proof.~} Since $(\mathcal {A}, [\cdot,\cdot]_\a, \a)$
is a Hom-Lie algebra and $(\mathcal {A}, \circ_\a, \a)$ is a
Hom-Novikov algebra, it remains to check the compatibility condition
\eqref{HGB}. As $\a$ is an algebra homomorphism of $(\mathcal {A}, [\cdot,\cdot],\circ, \a)$ and by \eqref{GB},
we have
\begin{eqnarray*} &&{[x
\circ_\a y, \a(z)]_\a -[x\circ_\a z, \a(y)]_\a +[x,y]_\a\circ_\a
\a(z)-[x,z]_\a \circ_\a \a(y)- \a(x)\circ_\a [y,z]_\a}\\
&&\ \ \ \ = \a^2\big([x \circ y, z] -[x\circ z, y] +[x,y]\circ
z-[x,z] \circ y- x\circ [y,z]\big)=0,
\end{eqnarray*}
which concludes the proof.
 \QED\vskip8pt
\begin{exam} \label{ex3}{\rm Suppose that $\D$ is an additive subgroup of $\mathbb{C}$
and denote $\Gamma=\{0,1,2,\cdots\}$. Let $\A_{\D,\Gamma}$ be the vector space with a
basis $\{x_{a,j}|(a,j)\in\Delta\times\Gamma\}$. For any given $\xi\in\C$, define an algebra operation on
$\mathcal{A}_{\D,\Gamma}$ by
\begin{eqnarray*}
x_{a,i}\c x_{b,j}=(b+\xi)x_{a+b,i+j}+jx_{a+b,i+j-1},\ \
\mbox{for} \ a,b\in\Delta, \ i,j\in\Gamma.
\end{eqnarray*}
It was proved in \cite{Xu2} that $(\A_{\D,\Gamma},\c)$ forms a
simple Novikov algebra. According to Gel'fand and Dorfman's
statement in \cite{GD}, there is a Gel'fand-Dorfman
bialgebra $(\A_{\D,\Gamma},[\cdot,\d]^-, \c)$ with
\begin{eqnarray*}
{[x_{a,i},x_{b,j}]^-}&=& x_{a,i}\c x_{b,j}-x_{b,j}\c  x_{a,i}\\
&=&{(b-a)x_{a+b,i+j}+(j-i)x_{a+b,i+j-1}},
\end{eqnarray*}
for $a,b\in\D$, $i,j\in\Gamma.$ Define
a linear endomorphism $\alpha$ of $\mathcal{A}_{\D,\Gamma}$ by
\begin{eqnarray}
\alpha(x_{a,j})=e^{a}x_{a,j}, \ \mbox{for}\ a\in\Delta.
\end{eqnarray}
It is easy to check that $\alpha$ is an algebra
homomorphism of $(\A_{\D,\Gamma},[\cdot,\d]^-, \c)$. By Theorem
\ref{p-2}, we get a Hom-GD bialgebra
$(\A_{\D,\Gamma},[\cdot,\d]^-_{\a}, \c_\a,\a)$. }
\end{exam}

 Let $(\A,\cdot)$ be a commutative associative algebra with a derivation $D$. It was proved that the new operation
 \begin{eqnarray}\label{Gd1}
x\circ y=x\d D(y)+\lambda x\d y, \ \mbox{for}\
x,y\in\mathcal {A},
\end{eqnarray}
equips $\mathcal {A}$ with a structure of Novikov algebra
for $\lambda=0$ by Gel'fand and Dorfman \cite{GD}, for
$\lambda\in\mathbb{F}$ by Filipov \cite{F}, and for a fixed element
$\lambda\in\mathcal {A}$ by Xu \cite{Xu1}. Yau \cite{Yau} gave an analogous construction of Hom-Novikov algebras (in $\l=0$ setting) from Hom-associative algebras introduced in \cite{MS}.

\begin{defi} {\rm A  Hom-associative algebra is a vector space $V$ with a bilinear map
$\mu: V \times V \longrightarrow V$ and a linear map $\alpha: V
\longrightarrow V$, such that
\begin{eqnarray}\label{hom-asso}
\mu(\alpha(x),\mu(y,z)) =
\mu(\mu(x,y),\alpha(z)),\ \mbox{for}\ x, y, z\in\mathcal{A}. \end{eqnarray}}
\end{defi}

A Hom-associative algebra
$(\mathcal {A},\mu,\alpha)$ is called {\it commutative }if
$\mu(x,y)=\mu(y,x)$ for $x,y\in\mathcal{A}$. A {\it derivation} on a Hom-associative algebra is defined in the usual way.

The following result presents our third construction of Hom-GD bialgebras
from a commutative Hom-associative algebra with a derivation.
\begin{theo}\label{twist}
Let $(\mathcal {A},\d, \alpha)$ be a commutative Hom-associative
algebra with a derivation $D$ such that $\alpha D=D\a$. Then
$(\mathcal {A},[\d,\d]^-, \circ, \alpha)$ forms a Hom-GD algebra,
where $\circ$ is defined by (\ref{Gd1}) and
\begin{eqnarray}\label{def1}
[x,y]^-=x\d D(y)-y\d D(x),\ \mbox{for}\ x, y \in\mathcal {A}.\end{eqnarray}
\end{theo}
\noindent{\it Proof.~} By \eqref{Gd1}, \eqref{hom-asso} and the fact that $\alpha D=D\a$, we have
\begin{eqnarray}
(x\c y)\c\a(z)&=&(x\d D(y)+\lambda x\d y)\c\a(z)\nonumber\\
&=&(x\d D(y)+\lambda x\d y)\d D(\a(z))+\lambda (x\d D(y)+\lambda x\d
y)\d \a(z)\nonumber\\&=& \a(x)\d (D(y)\d D(z))+\lambda \a(x)\d
\big(y\d D(z)+D(y)\d z\big)+\lambda^2 \a(x)\d (y\d z). \ \ \ \ \label{1}
\end{eqnarray}
Similarly, we have
\begin{eqnarray}\label{2}
(x\c z)\c\a(y)= \a(x)\d (D(z)\d D(y))+\lambda\a(x)\d\big(z\d
D(y)+D(z)\d y\big)+\lambda^2 \a(x)\d (z\d y).
\end{eqnarray}
Combining \eqref{1} with \eqref{2} and using commutativity of
$(\A,\d,\a)$, we get
\begin{eqnarray*}
(x\c y)\c\a(z)=(x\c z)\c\a(y),
\end{eqnarray*}
which proves \eqref{Hom-Nov2+}. By commutativity and Hom-associativity of $(\A,\d,\a)$, one can obtain
\begin{eqnarray}\label{aaaa}
(x\c y)\c\a(z)-\a(x)\c(y\c z)
=-(x\d y) \d \a(D^2(z))-\lambda(x\d y) \d \a(D(z)),
\end{eqnarray}
which implies  \eqref{Hom-Nov1+}, because the right-hand side of \eqref{aaaa} is symmetric in $x$ and $y$. So $(\mathcal {A}, \circ, \alpha)$
is a Hom-Novikov algebra. By \eqref{def1} and commutativity of $(\A,\d)$, we have $[x,y]^-=x\c y-y\c x$. It follows from Theorem \ref{th1} that $(\mathcal {A},[\d,\d]^-, \circ, \alpha)$
is a Hom-GD algebra.\QED\vskip8pt

From the proof of the above theorem, we have the following result, generalizing Yau's construction of Hom-Novikov algebras (see \cite[Theorem 1.2]{Yau}).
\begin{coro} Let $(\mathcal {A},\d, \alpha)$ be a commutative Hom-associative algebra
with a derivation $D$ such that $\alpha D=D\a$. For any fixed number $\lambda$, define
 \begin{eqnarray}
x\circ y=x\d D(y)+\lambda x\d y, \ \mbox{for}\
x,y\in\mathcal {A}.
\end{eqnarray}
Then $(\mathcal {A},\circ, \alpha)$ is a Hom-Novikov algebra.
\end{coro}

Moreover, we have
\begin{coro}\label{coro} Let $(\A,\d)$ be a commutative associative algebra with
an algebra endomorphism $\a$ and
a derivation $D$, such that $D\a=\a D$. For any fixed number $\lambda$, define
\begin{eqnarray}\label{gd3}
x\circ y=\a(x\d D(y))+\lambda\a(x\d y), \hspace{0.2cm}\mbox{for}\
x,y\in\mathcal {A}.
\end{eqnarray}
Then $(\A,[\d,\d]^-, \c,\a)$
forms a Hom-GD bialgebra, where
\begin{eqnarray}
[x,y]^-=\a(x\d D(y))-\a(y\d D(x)), \ \mbox{for}\
x,y\in\mathcal {A}.
\end{eqnarray}
\end{coro}
\noindent{\it Proof.~} For clarity, write $x\d_\a y=\a(x\d y)$ for
$x,y\in\A$. Because $\a$ is an algebra endomorphism of $(\A,\d)$, we have
$(\A,\d_\a,\a)$ is a commutative Hom-associative algebra. Since
$\a D=D\a$, $D$ is also a derivation of $(\A,\d_\a)$. We can rewrite
\eqref{gd3} as
\begin{eqnarray*}
x\circ y=x\d_\a D(y)+\lambda x\d_\a y \hspace{0.25cm}\mbox{for}\
x,y\in\mathcal {A}.
\end{eqnarray*}
According to Theorem \ref{twist}, $(\A,[\d,\d]^-, \c,\a)$ forms a
Hom-GD bialgebra.\QED\vskip8pt

The following three examples are applications of Corollary
\ref{coro}.
\begin{exam} {\rm Consider the polynomial algebra $\mathbb{C}[x]$. Let $D=\frac{\rm d}{{\rm d }x}$ be the differential operator and $\a$ an endomorphism of $\mathbb{C}[x]$ determined by
\begin{eqnarray*}
\a(x^n)=(x+c)^n , \ \mbox{with}\ c\in\mathbb{C}.
\end{eqnarray*}
One has $\a D=D\a$ (cf.
\cite{Yau}). By Corollary
\ref{coro}, we have a Hom-GD bialgebra $(\mathbb{C}[x], [\cdot, \cdot]^-,\circ,\a)$ with
\begin{eqnarray*}
f(x)\circ g(x)&=&f(x+c)D(g(x+c))+\l f(x+c)g(x+c),\\ {[f(x),g(x)]^-}&=&f(x+c)D(g(x+c))-g(x+c)D(f(x+c)),
\end{eqnarray*}
for all $f(x),g(x)\in\mathbb{C}[x]$, where $\lambda$ is a fixed number. }\end{exam}
\begin{exam}{\rm Let $A$ be a commutative associative algebra with a nilpotent derivation
$D$, namely, there exists a
positive integer $n$ such that $D^n=0.$ The formal exponential
map
\begin{eqnarray*}
\a=e^D={\rm id}+D+\frac1{2!}D^2+\cdots+\frac1{(n-1)!}D^{n-1}
\end{eqnarray*}
is an algebra automorphism of $A$, and $\a D=D\a$ (cf.
\cite{Yau}). For a fixed
number $\lambda$, the operations
\begin{eqnarray*}
x\circ y=\a(x\d D(y))+\lambda\a(x\d y),\ {[x,y]^-}={\a(x\d
D(y))-\a(y\d D(x))}, \ \mbox{for}\  x,y\in\mathcal{A},
\end{eqnarray*}
give a Hom-GD bialgebra $(A,[\cdot, \cdot]^-,\circ,\a)$ by
Corollary \ref{coro}. }
\end{exam}

\begin{exam}\label{ex2}{\rm Let $\Delta$ be a nonzero abelian subgroup of
$\mathbb{C}$. Suppose that $f$ is a nontrivial homomorphism of
$\Delta$ into the additive group of $\mathbb{C}$. Consider the vector
space $\mathcal{A}$ with a basis $\{x_a|a\in\Delta\}$. Define an
operation on $\mathcal{A}$ by
\begin{eqnarray*}
x_a\cdot x_b=x_{a+b}, \hspace{0.3cm}\mbox{for all} \ a,b\in\Delta.
\end{eqnarray*}
Then $(\mathcal{A}, \cdot)$ is a commutative associative  algebra.
Define a linear map $\alpha:\mathcal{A}\longrightarrow\mathcal{A}$
by
\begin{eqnarray*}
\alpha(x_a)=e^{a}x_{a}, \hspace{0.3cm}\mbox{for all} \ a\in\Delta.
\end{eqnarray*}
Then $\alpha$ is an algebra homomorphism of $(\mathcal{A},\cdot)$.
Define another linear map
$\partial:\mathcal{A}\longrightarrow\mathcal{A}$ by
\begin{eqnarray*}
\partial(x_a)=f(a)x_{a}, \hspace{0.3cm}\mbox{for all} \ a\in\Delta.
\end{eqnarray*}
It is easy to check that $\partial$ is a derivation of
$(\mathcal{A}, \cdot)$ and $\a\partial=\partial\a $. By Corollary
\ref{coro}, we obtain a Hom-GD bialgebra $(A,[\cdot,
\cdot]^-,\circ,\a)$ with
\begin{eqnarray*}
x_a\circ x_b=(f(b)+\lambda)e^{a+b}x_{a+b},\ \ {[x_a,x_b]^-}=f(a+b)e^{a+b}x_{a+b},
\end{eqnarray*}
for all $a,b\in\Delta$, where $\lambda$ is a fixed number.}
\end{exam}

Our fourth construction of Hom-GD bialgebras is related to
 the notion of Hom-Poisson algebra, which was introduced
in the study of deformation theory of Hom-Lie algebras \cite{MS1}.
\begin{defi}{\rm A Hom-Poisson algebra is a vector space $\A$ equipped with
two operations $\cdot$ and $[\cdot,\cdot]$, and a linear
endomorphism $\a$, such that $(\A,\cdot,\a)$ is a commutative Hom-associative algebra,
$(\A,[\cdot,\cdot],\a)$ is a Hom-Lie algebra, and the following relation holds for $x, y, z\in \A$:
\begin{eqnarray}\label{p}
[\a(x),y\d z] = \a(y)\d[x, z] + \a(z)\d[x, y].
\end{eqnarray}}
\end{defi}

Notice that \eqref{p} can be equivalently reformulated as
\begin{eqnarray}\label{p-1}
[x\d y,\a(z)] = [x, z]\d\a(y)+ \a(x)\d[y, z].
\end{eqnarray}
By setting $\a={\rm id}$ in above definition, we recover Lie-Poisson algebras, which appear naturally in Hamiltonian
mechanics, and are central in the study of quantum groups.

 The following result can
be considered as a Hom-version of [\ref{X1}, Theorem 3.2].

\begin{theo} \label{th2} Let $(\A, \d, [\cdot,\cdot],\a)$ be a Hom-Poisson
algebra and $D$ a derivation of $(\A, \d)$, such that $D\a=\a D$
and
\begin{eqnarray}\label{gd2}
D([x,y])=[D(x),y]+[x,D(y)]+\lambda[x,y],\ \mbox{for} \ x,y\in\A,
\end{eqnarray}
where $\lambda$ is a fixed number. Define a new
operation $\c$ by
\begin{eqnarray}\label{hhh} x\circ y=x\d
D(y)+\lambda x\d y, \ \mbox{for}\ x,y\in\mathcal {A}.
\end{eqnarray} Then $(\mathcal {A},[\d,\d],\circ,\a)$ forms a Hom-GD
bialgebra.
\end{theo}
\noindent{\it Proof.~} By Theorem \ref{twist}, $(\mathcal
{A},\circ,\a)$ is a Hom-Novikov algebra. By \eqref{p}--\eqref{hhh} and the fact that $D$ is a derivation of
$(\A,\d)$ commuting with $\a$, we get
\begin{eqnarray*}
&&{[x \circ y, \a(z)] -[x\circ z, \a(y)] +[x,y]\circ
\a(z)-[x,z] \circ \a(y)- \a(x)\circ [y,z]}\\
&&\ \ \ = [x\d D(y)+\lambda x\d y, \a(z)]-[x\d D(z)+\lambda x\d
z,\a(y)]+[x,y]\d \a(D(z))+\lambda [x,y]\d\a(z)\\&&\ \ \ \ \ \ -[x,z]\d
\a(D(y))-\lambda [x,z]\d\a(y)-\a(x)\d D([y,z])-\lambda\a(x)\d [y,z]
\\&&\ \ \ =[x\d D(y), \a(z)]-[x\d D(z)+\lambda x\d
z,\a(y)]+[x,y]\d \a(D(z))+\lambda [x,y]\d\a(z)\\&&\ \ \ \ \ \ -[x,z]\d
\a(D(y))-\a(x)\d \big([D(y),z]+[y,D(z)]+\lambda[y,z]\big)\\&&\ \ \ =
\big([x\d
D(y),\a(z)]-[x,z]\d\a(D(y))-\a(x)\d[D(y),z]\big)\\&&\ \ \ \ \ \ -\big([x\d
D(z),\a(y)]-[x,y]\d\a(D(z))-\a(x)\d[D(z),y]\big)= 0,
\end{eqnarray*}
which proves \eqref{HGB} and the theorem. \QED\vskip8pt
As it was explained in \cite{X1}, the above construction is related to Lie superalgebras of
Hamiltonian type and Contact type \cite{Ka2}.

\vs{10pt}\

\cl{\bf\S4. \ Hom-Lie conformal algebras
}\setcounter{section}{4}\setcounter{equation}{0}\setcounter{theo}{0}
\vskip8pt

In this section, we introduce the notion of Hom-Lie conformal
algebra and give a general construction of Hom-Lie conformal
algebras from formal distribution Hom-Lie algebras, analogous to the construction of Lie conformal algebras from formal distribution Lie algebras.

\begin{defi}\label{HLCF}{\rm
A Hom-Lie conformal algebra is a $\C[\partial ]$-module $\mathcal
{R}$ equipped with a linear endomorphism $\a$ such that $\a\p=\p\a$, and a $\lambda$-bracket $[\cdot_\lambda\cdot]$ which defines a $\C$-bilinear map from $\mathcal {R}\otimes \mathcal {R}$
to $\mathcal {R}[\lambda]=\C[\l]\otimes \mathcal{R}$ such that the following axioms hold for $a, b, c\in \mathcal {R}$:
\begin{eqnarray}
[\partial a_\lambda b]&=&-\lambda[a_\lambda b],\ \ \ \ \ \mbox{(conformal\  sesquilinearity)}\label{HL1}\\
{[a_\lambda b]} &=& -[b_{-\lambda-\partial}a],\ \ \ \ \ \ \mbox{(skew-symmetry)}\label{HL2}\\
{[\a(a)_\lambda[b_\mu c]]}&=&[[a_\lambda b]_{\lambda+\mu
}\a(c)]+[\a(b)_\mu[a_\lambda c]].\ \ \mbox{(Hom-Jacobi \
identity)}\label{HL3}
\end{eqnarray}
}\end{defi}
\begin{rema}{\rm We recover Lie conformal algebras when $\a={\rm id}.$ In addition, we can associate to any Lie conformal algebra a Hom-Lie conformal algebra structure by
taking $\a=0$. Such a Hom-Lie conformal algebra is called {\it trivial}.}\end{rema}
\begin{rema}{\rm
It follows from \eqref{HL1} and \eqref{HL2} that
\begin{eqnarray}\label{HL4}
[a_{\l}\p b]=(\p+\l)[a_{\l} b],  \ \ \mbox{for} \
a,b\in\mathcal{R}.
\end{eqnarray}
Thus $\p$ acts as a derivation on the $\l$-bracket. }\end{rema}

A Hom-Lie
conformal algebra $\mathcal {R}$ is called {\it finite} if
$\mathcal {R}$ is a finitely generated $\C[\partial]$-module. The
{\it rank} of $\mathcal {R}$ is its rank as a $\C[\partial]$-module. Let
$m$ be a positive integer. $\mathcal {R}$ is called a {\it Hom-Lie conformal algebra of degree $m$} if for any $a,
b \in \mathcal {R}$, there exist $w_{i,j}\in\mathcal{R}$ such that
\begin{eqnarray}
[a_\l b]=\mbox{$\sum\limits_{j\geq0,i+j< m}$}\p^i w_{i,j}\l^{j}.
\end{eqnarray}
 In particular,  for $m=2$, there exist
$c_{0,0}, c_{1,0}, c_{0,1}\in\mathcal {R}$ such that
 \begin{eqnarray}\label{quadratic}
[a_\l b]=c_{0,0}+\p c_{1,0} +\l c_{0,1}.
\end{eqnarray}

Let $V$ be any vector space. Following \cite{ka}, a {\it $V$-valued formal distribution} is a formal series
\begin{eqnarray}a(z) =\mbox{$\sum\limits_{n\in\Z}$}a_{(n)}z^{-n-1},\end{eqnarray}
where $a_{(n)}={\rm Res}_z z^na(z)\in V$ is called the {\it
Fourier coefficient} of $a(z)$. The notion of {\it residue} is
taken formally in analogy with Complex Analysis.
The vector space of these formal distributions in one
variable is denoted by $V[[z,z^{-1}]]$. Similarly, a formal
distribution $a(z,w)$ in two variables is defined as a series $\sum_{m,n\in\Z}a_{m,n}z^{-m-1}w^{-n-1}$ and the space of these
series is denoted by $V[[z,z^{-1},w,w^{-1}]]$. The {\it delta
distribution} is the $\C$-valued formal distribution
\begin{eqnarray}\label{deta}
\delta(z,w)=\mbox{$\sum\limits_{n\in\Z}$}z^{-n-1}w^n.
\end{eqnarray}
It enjoys the following properties:
\begin{eqnarray}
\partial_z\delta(z,w)=-\partial_w\delta(z,w),\ \
(z-w)^m\partial_w^n\delta(z,w)=0 \ \ \mbox{for}\
m>n.\label{deta33}
\end{eqnarray}

It is known that operator product expansions (OPE's) are widely
used in conformal field theory. The fundamental notion beneath
it is the locality of formal distributions in two variables.
\begin{defi}\rm Let $a(z,w)\in V[[z,z^{-1}, w,w^{-1}]].$ The formal
distribution $a(z,w)$ is called local if there exists a positive
integer $N$ such that $(z-w)^Na(z,w)=0.$
\end{defi}

The following decomposition theorem is due to Kac \cite{ka}.
\begin{theo} \label{decom} Let $a(z,w)\in V[[z,z^{-1},w,w^{-1}]]$ be a local formal distribution. Then $a(z,w)$ can be written as a finite sum of $\delta(z,w)$ and its derivatives:
\begin{eqnarray}
a(z,w)=\mbox{$\sum\limits_{j\in\Z^+}$}c^j(w)\frac{\p_w^{j}\delta(z,w)}{j!},
\end{eqnarray}
where
$c^j(w)={\rm Res}_z(z-w)^ja(z,w)\in V[[w,w^{-1}]].$
In addition, the converse is true.
\end{theo}

Denote
\begin{eqnarray}
e^{\lambda(z-w)}= \mbox{$\sum\limits_{j\in \Z^+
}$}\frac{\lambda^j}{j!}(z-w)^j\in
\mathbb{C}[[z,z^{-1},w,w^{-1}]][[\lambda]].
\end{eqnarray}
The {\it formal
Fourier transform} of $a(z,w)\in
V[[z,z^{-1},w,w^{-1}]]$ is defined
by
\begin{eqnarray}\label{fourier}
F^\lambda_{z,w}a(z,w)={\rm Res}_ze^{\lambda(z-w)}a(z,w).
\end{eqnarray}
Notice that $F^\lambda_{z,w}$ is a linear map from $
V[[z,z^{-1},w,w^{-1}]]$ to $V[[w,w^{-1}]][[\lambda]].$ Indeed, we have
\begin{eqnarray}\label{fourier2}
F^\lambda_{z,w}a(z,w)=\mbox{$\sum\limits_{j\in\mathbb{Z}^+}$}\frac{\lambda^j}{j!}c^j(w),\
\ \ \mbox{where}\ c^j(w)={\rm Res}_z(z-w)^ja(z,w).
\end{eqnarray}
If $a(z,w)$ is local, then $F^\lambda_{z,w}a(z,w)\in
V[[w,w^{-1}]][\lambda]$, and
\begin{eqnarray}\label{fourier3}
F^\lambda_{z,w}\p_z a(z,w)=-\l
F^\lambda_{z,w}a(z,w)=[\p_w,F^{\l}_{z,w}], \ \ \
F^\lambda_{z,w}a(w,z)=F^\mu_{z,w}a(z,w)|_{\mu=-\l-\p_w}.
\end{eqnarray}

Let $(\LL,[\d,\d],\a)$ be a Hom-Lie algebra.
We can extend the Hom-Lie bracket to the commutator
between two $\LL$-valued formal distributions $a(z)=\mbox{$\sum_{m\in\Z}$}a_mz^{-m-1}$ and $b(w)=\sum_{n\in\Z}b_nw^{-n-1}$
by
\begin{eqnarray}\label{LB2}
[a(z),b(w)]=
\mbox{$\sum\limits_{m,n\in\Z}$}[a_{(m)},b_{(n)}]z^{-m-1}w^{-n-1}\in
\LL[[z,z^{-1},w,w^{-1}]].
\end{eqnarray}
If $[a(z),b(w)]$ is local, we say $a(z)$ and $b(w)$ are {\it mutually local}, or $\big(a(z),b(w)\big)$ is a {\it local pair}.
Taking the formal Fourier transform  $F^\lambda_{z,w}$ defined by \eqref{fourier} on
both parts of \eqref{LB2}, we can define a $\mathbb{C}$-bilinear map $[\d\,_{\lambda}\d]$, called a {\it $\l$-bracket}, from  $\LL[[w,w^{-1}]]\otimes
\LL[[w,w^{-1}]]$ to $\LL[[w,w^{-1}]][[\lambda]]$ by
\begin{eqnarray}\label{lamda}
[a(w)_\lambda b(w)]= F^\lambda_{z,w}[a(z),b(w)].
\end{eqnarray}
In the sequel, the $\lambda$-bracket between $a(w)$ and $b(w)$ will be simply denoted by
$[a_\lambda b]$.

Let $\big(a(z),b(w)\big)$ be a local pair of $\LL$-valued formal distributions. By Theorem \ref{decom}, we have
\begin{eqnarray}\label{sum+}
[a(z),b(w)]=\mbox{$\sum\limits_{j\in\Z^+}$}a(w)_{(j)}b(w)\frac{\p_w^{j}\delta(z,w)}{j!},
\end{eqnarray}
where
\begin{eqnarray}\label{jsum+}
a(w)_{(j)}b(w)={\rm Res}_z(z-w)^j[a(z),b(w)]
\end{eqnarray}
is called the {\it $j$-product} of $a(w)$ and $b(w)$. This product
will be simply denoted by $a_{(j)}b$. By \eqref{lamda} and \eqref{sum+}, the $\lambda$-bracket
is related to the $j$-products as follows:
\begin{eqnarray}\label{lj}
[a_\lambda
b]=\mbox{$\sum\limits_{j\in\Z^+}$}\frac{\lambda^j}{j!}(a_{(j)}b).
\end{eqnarray}

By defining the action of $\p$ on
$a(z)=\mbox{$\sum_{m\in\Z}$}a_{(m)}z^{-m-1}\in\LL[[z,z^{-1}]]$ by $(\p a)(z)=\p_z(a(z))$, and the action
of $\a$ on $a(z)$ by $\a(a)(z)=\mbox{$\sum_{m\in\Z}$}\a(a_{(m)})z^{-m-1}$, we have
the following results.
\begin{prop} \label{p*} The $\l$-bracket satisfies the following properties:
\begin{itemize}
\item[\rm (1)]$[\partial a_\lambda
b]=-\lambda[a_\lambda b],\ \ {[a_\lambda
\partial b]}=(\p+\lambda)[a_\lambda b].$
\item[\rm(2)] if $(a(z),b(w))$ is a local
pair, then ${[a_\lambda b]}=-[b_{-\lambda-\partial}a].$
\item[\rm(3)] ${[\a(a)_\lambda[b_\mu c]]}=[[a_\lambda b]_{\lambda+\mu
}\a(c)]+[\a(b)_\mu[a_\lambda c]].$
\end{itemize}
\end{prop}
\noindent{\it Proof.~} (1) and (2) follow from \cite [Proposition 2.3 (a), (b), (c)]{ka}.
By \eqref{JJ} and \eqref{LB2}, we have
\begin{eqnarray}\label{222}
[\a(a)(z),[b(x),c(w)]]=[\a(b)(x),[a(z),c(w)]]+[[a(z),b(x)],\a(c)(w)].
\end{eqnarray}
This together with \eqref{lamda} gives
\begin{eqnarray*}
[\a(a)_\lambda[b_\mu
c]]&=&F^{\l}_{z,w}[\a(a)(z),F^{\mu}_{x,w}[b(x),c(w)]]=F^{\l}_{z,w}F^{\mu}_{x,w}[\a(a)(z),[b(x),c(w)]]\\
&=&F^{\l}_{z,w}F^{\mu}_{x,w}\big([\a(b)(x),[a(z),c(w)]]+[[a(z),b(x)],\a(c)(w)]\big).
\end{eqnarray*}
Because of
$[[a(z),b(x)],\a(c(w))]\in\LL[[z,z^{-1},x,x^{-1},w,w^{-1}]]$, we can
use the property of the formal Fourier transform, according to which
$F^{\l}_{z,w}F^{\mu}_{x,w}=F^{\l+\mu}_{x,w}F^{\l}_{z,x}$. Hence,
\begin{eqnarray*}
[\a(a)_\lambda[b_\mu
c]]&=&F^{\l}_{z,w}F^{\mu}_{x,w}\big([\a(b)(x),[a(z),c(w)]]\big)+F^{\l+\mu}_{x,w}F^{\l}_{z,x}\big([[a(z),b(x)],\a(c)(w)]\big)\\
&=&F^{\mu}_{x,w}\big([\a(b)(x),F^{\l}_{z,w}[a(z),c(w)]]\big)+F^{\l+\mu}_{x,w}\big([F^{\l}_{z,x}[a(z),b(x)],\a(c)(w)]\big)\\
&=&[\a(b)_\mu[a_\lambda c]]+[[a_\lambda b]_{\lambda+\mu }\a(c)],
\end{eqnarray*}
which proves (3).\QED\vskip8pt

By \eqref{lj}, Proposition \ref{p*} is translated in terms
of $j$-products as follows.
\begin{prop} \label{p**}The $j$-products satisfy the following equalities.
\begin{itemize}
\item[\rm(1)]${\partial a_{(n)} b}=-na_{(n-1)} b,\ {a_{(n)} \partial b}=\p(a_{(n)}b)+n
a_{(n-1)}b$.
\item[\rm(2)]${{a_{(n)}b}} = -\mbox{$\sum_{i\geq 0}$}(-1)^{n+i}\frac{1}{i!}\p^i
b_{(n+i)}a,$ if $(a(z),b(w))$ is a local pair.
\item[\rm(3)]${\a(a)_{(m)}b_{(n)} c}={\a(b)_{(n)}a_{(m)}
c}+\mbox{$\sum\limits_{i=0}^{m}$}\mbox{${m\choose
i}$}(a_{(i)}b)_{(m+n-i)}\a(c).$
\end{itemize}
\end{prop}

A subset $F\subset \LL[[z,z^{-1}]]$ is called a {\it local family}
of $\LL$-valued formal distributions if all pairs of its
constituents are local. The following notion of formal
distribution Hom-Lie algebra can be seen as a Hom-analogue of the notion of
formal distribution Lie algebra introduced by Kac \cite{ka}.

\begin{defi}{\rm Let $(\LL,[\d,\d],\a)$ be a Hom-Lie algebra.
If there exists a local family $F$ of $\LL$-valued formal
distributions, with their Fourier coefficients generating the whole
$\LL$, then $F$ is said to endow $\LL$ with a structure of
formal distribution Hom-Lie algebra. In this case, we denote
$\LL$ by $(\LL,F)$ to emphasize the role of $F$.}
\end{defi}

Here is a simple example of formal distribution Hom-Lie algebras.
\begin{exam}\label{ex4}
{\rm  Let $(\LL, [\d,\d],\a)$ be a Hom-Lie algebra. Denote
$$F=\{g(z)=\mbox{$\sum\limits_{n\in\Z}$}g z^{-n-1}|g\in \LL\}.$$
For any $g,h\in \LL$, we have
\begin{eqnarray*}
[g(z),h(w)]=\mbox{$\sum_{m,n}$}[g,h]z^{-m-1}w^{-n-1}
=[g,h]\mbox{$\sum_{m,k}$}z^{-m-1}w^{-k+m-1}
=[g,h](w)\delta(z,w).
\end{eqnarray*}
Thus $(g(z),h(w))$ is a local pair, and $(\LL, F)$ is a formal distribution
Hom-Lie algebra.
 }\end{exam}

The previous example shows that we can trivially associate to any
Hom-Lie algebra $\LL$ a structure of formal distribution Hom-Lie
algebra. In the following, we shall put other restrictions to
the local family $F$, such that we can endow the formal distribution Hom-Lie algebra $(\LL,
F)$ with a Hom-Lie conformal
algebra structure.

\begin{defi} {\rm Let $(\LL, [\d,\d],\a)$ be a Hom-Lie algebra. A local family
$F\subset \LL[[z,z^{-1}]]$ is called a
conformal family if it is closed under their $j$-products and invariant under the actions of $\p$ and $\a$.}
\end{defi}

It is known from \cite{ka} that if $(\LL, F)$ is a formal
distribution Lie (super)algebra, one can always include $F$ in the
minimal conformal family $\bar F$, such that $(\LL, F)$ can be given a Lie conformal algebra structure $\mathcal{R}=\C[\p]\bar F,$ with
$[a_\l b]=F^{\l}_{z,w}[a(z),b(w)]$, $\p=\p_z$. For a formal
distribution Hom-Lie algebra $(\LL, F)$, we also have a Hom-Lie conformal
algebra $\mathcal{R}=\C[\p]\bar F$ by Propositions \ref{p*} and \ref{p**}. Indeed, it suffices to extend $\a$ to $\mathcal{R}=\C[\p]\bar F$ by $\a(f(\p)u)=f(\p)\a(u)$, for $f(\p)\in\C[\p],$ $u\in \bar F$.

The following are examples of Hom-Lie conformal algebras.
\begin{exam}{\rm  Let $(\LL, [\d,\d],\a)$ be a Hom-Lie algebra. Denote
by $\hat \LL=\LL\otimes \C[t,t^{-1}]$ the affization of $\LL$ with
\begin{eqnarray*}
[u\otimes t^m, v\otimes t^n]=[u,v]\otimes t^{m+n},\ \ \mbox{for} \
u, v\in L,\ m,n\in\Z.
\end{eqnarray*}
Extend $\a$ to $\hat \LL$ by $\a(u\otimes t^m)=a(u)\otimes t^m$.
Then $(\hat \LL, [\d,\d],\a)$ is a Hom-Lie algebra. Let
$$u(z)=\mbox{$\sum\limits_{n\in\Z}$}(u\otimes t^n)z^{-n-1}.$$
We have
\begin{eqnarray*}
[u(z),v(w)]&=&\mbox{$\sum\limits_{m,n}$}([u,v]\otimes t^{m+n})z^{-m-1}w^{-n-1}\\
&=&\mbox{$\sum\limits_{m,k}$}([u,v]\otimes t^k)z^{-m-1}w^{-k+m-1}\\
&=&[u,v](w)\delta(z,w),
\end{eqnarray*}
which is equivalent to $u_{(0)}v=[u,v]\in \LL$, and $u_{(n)}v=0$ for $n>0$, namely,
$\LL$ is closed under all the $j$-products.
Thus, we get a Hom-Lie conformal algebra $\mathcal{R}=\C[\p]\LL$
with
\begin{eqnarray}\label{1111}
[u_{\l}v]=[u,v], \ \mbox{for} \ u,v\in \LL.
\end{eqnarray}
Indeed, we can extend the $\lambda$-bracket to the whole $\mathcal{R}$
by
\begin{eqnarray}\label{lamda22+2}
[f(\partial) u_\lambda h(\partial)
v]= f(-\lambda)h(\partial+\lambda)[u_\lambda v],
\end{eqnarray}
and extend $\a$ to a linear map of $\mathcal{R}$ by
\begin{eqnarray}\label{lamda22+2a}\a(f(\p)u)=f(\p)\a(u),\end{eqnarray} for any $f(\partial),h(\partial)\in \mathbb{C}[\p]$,
$u,v\in\mathcal{R}$.
The conformal sesquilinearity \eqref{HL1} is naturally satisfied by
(\ref{lamda22+2}), and (\ref{lamda22+2a}) gives $\a\p=\p\a$. Thus it suffices to check \eqref{HL2} and \eqref{HL3} on the generators.
Since the $\lambda$-bracket \eqref{1111} is defined by the Hom-Lie bracket on $\LL$, axioms \eqref{HL2} and \eqref{HL3} hold. So $\mathcal{R}=\C[\p]\LL$ is a Hom-Lie conformal algebra, which is viewed as a {\it current-like  Hom-Lie conformal algebra}. }\end{exam}

\begin{rema}\label{remark}{\rm If $\mathcal{R}$ is a free $\C[\p]$-module with the $\l$-bracket
 $[\d_{\l}\d]$ defined on a $\C[\p]$-basis of $\mathcal{R}$ such that \eqref{HL2}
and \eqref{HL3} hold, there is a unique extension of this
$\l$-bracket via \eqref{HL1} and \eqref{HL4} to the whole
$\mathcal{R}$ (as shown in \eqref{lamda22+2}). It is easy to show that
\eqref{HL2} and \eqref{HL3} also hold for this extension. Thus,
$\mathcal{R}$ is equipped with a Hom-Lie conformal algebra
structure. In the sequel, we shall often describe Hom-Lie conformal
algebra structures on free $\C[\p]$-modules by giving the
$\l$-bracket on a fixed $\C[\p]$-basis. }\end{rema}

\begin{exam}\label{exm2}{\rm Recall that the Virasoro conformal algebra is a free $\C[\p]$-module $Vir=\C[\p]L$
generated by one symbol $L$ such that
 \begin{eqnarray}
 [L_\l L]=(\p+2\l)L.
 \end{eqnarray}
Define $\a(L)=f(\p)L$, with $0\neq f(\p)\in\C[\p]$. If $(Vir, [\cdot_\l\cdot],\a)$ forms a Hom-Lie conformal algebra, we have
\begin{eqnarray*}
[\a(L)_\l[L_\mu L]]=[\a(L)_\mu[L_\l L]]+[[L_\l L]_{\l+\mu}\a(L)],
\end{eqnarray*}
which is equivalent to the following equality:
\begin{eqnarray}\label{3+}
&&f(-\l)(\p+\l+2\mu)(\p+2\l)\nonumber\\
&&\ \ \ \ \ \ \ \ \ \, =f(-\mu)(\p+\mu+2\l)(\p+2\mu)+(\l-\mu)f(\l+\mu+\p)(\p+2\l+2\mu).
\end{eqnarray}
By \eqref{3+}, the highest degree of $\p$ in $f(\p)$ is at most $1$. We can write $f(\p)=a\p+b$, with $a,b\in\C$ and $(a,b)\neq (0,0)$. Comparing the coefficients of $\p^2$ in \eqref{3+}, we have
\begin{eqnarray}
a(\mu-\l)\p^2=a(\l-\mu)\p^2,
\end{eqnarray}
which implies $a=0$ and thus $f(\p)=b\neq 0$. Consequently, we obtain a Hom-Lie
conformal algebra $Vir=\C[\p]L$ with $[L_{\l}L]=(\p+2\l)L$, and
$\a(L)=b L$, where $b$ is a nonzero number. We call it the {\it
Virasoro-like Hom-Lie conformal algebra}.}\end{exam}

 With $$[a_{\lambda}b]=\mbox{$\sum\limits_{j\in\Z^+
}$}\frac{\lambda^j}{j!}(a_{(j)}b),$$ we can give a
 equivalent definition of a Hom-Lie conformal algebra.
 \begin{defi}{\rm A Hom-Lie conformal algebra is a $\C[\partial ]$-module
$\mathcal {R}$ equipped with a linear map $\a$ such that
$\a\p=\p\a$, and equipped with infinitely many $\C$-bilinear
$j$-products $(a,b)\rightarrow a_{(j)}b$ with $j\in\Z^+$, such that
the following axioms hold for $a,b,c \in\mathcal {R}$, $m,n
\in\Z^+$:
\begin{eqnarray}
{a_{(n)}b}&=&0 \ \ \mbox{for} \ n \ \mbox{sufficiently large},\\
{\partial a_{(n)} b}&=&-na_{(n-1)} b,\label{jHL1}\\
{{a_{(n)}b}} &=& -\mbox{$\sum_{i\geq 0}$}(-1)^{n+i}\frac{1}{i!}\p^i b_{(n+i)}a,\label{jHL2}\\
{\a(a)_{(m)}b_{(n)} c}&=&{\a(b)_{(n)}a_{(m)}
c}+\mbox{$\sum\limits_{i=0}^{m}$}\mbox{${m\choose
i}$}(a_{(i)}b)_{(m+n-i)}\a(c). \label{jjHL3}
\end{eqnarray}
}\end{defi}

\begin{rema}{\rm In terms of $j$-products, \eqref{HL4} is
equivalent to
\begin{eqnarray*}
{a_{(n)} \partial b}=\p(a_{(n)}b)+n a_{(n-1)}b,
\end{eqnarray*}
which together with \eqref{jHL1} shows that $\p$ acts as a derivation on the $j$-products.}
\end{rema}

\vs{6pt}\

\cl{\bf\S5. \ Equivalence
}\setcounter{section}{5}\setcounter{equation}{0}\setcounter{theo}{0}
\vskip8pt
In this section, we show that the affinization of a Hom-GD bialgebra forms a Hom-Lie algebra. Also, we prove that a Hom-GD bialgebra is equivalent to a Hom-Lie conformal algebra of degree $2$.

Let $\mathcal {A}$ be a vector space with two bilinear operations $[\d,\d]$ and
$\circ$, and a linear endomophism $\a$.
Denote
\begin{eqnarray*}
\mathcal {L}(\A)=\A\otimes\mathbb{C}[t,t^{-1}].
\end{eqnarray*}
 Define a bilinear operation $[-,-]$ on $\mathcal {L}(\A)$ 
by
\begin{eqnarray}\label{aff}
[u\otimes t^m, v\otimes t^n]=[u,v]\o t^{m+n}+m u\c v\o t^{m+n-1}-n
v\c u\o t^{m+n-1},
\end{eqnarray}
for all $u,v\in\A$, and $m,n\in\mathbb{Z}$. Moreover, define a linear map
$\varphi:\mathcal {L}(\A)\rightarrow \mathcal {L}(\A)$ by
\begin{eqnarray}\label{map}
\varphi(u\otimes t^m)=\a(u)\otimes t^m, \ \ \mbox{for}\ u\in \A, \ m\in\mathbb{Z}.
\end{eqnarray}

We have the following result.
\begin{theo}\label{h-h} $(\L,[-,-],\varphi)$
is a Hom-Lie algebra if and only if $(\mathcal
{A},[\d,\d],\circ,\a)$ is a Hom-GD bialgebra.
\end{theo}
\noindent{\it Proof.~} Assume that
$(\mathcal {A},[\d,\d],\circ,\a)$ is a Hom-GD bialgebra. By \eqref{aff}, the bracket $[-,-]$ on $\mathcal {L}(\A)$ is
skew-symmetric since the bracket $[\d,\d]$ on $\A$ is skew-symmetric.  For
$u,v,w\in\A$, $m,n,k\in\mathbb{Z}$, we have 
\begin{eqnarray}\label{jj}
&&[[u\otimes t^m, v\otimes t^n],\varphi(w\o t^k)]+[[v\o t^n,w\o
t^k],\varphi(u\o t^m)]+[[w\o t^k,u\o t^m],\varphi(v\o
t^n)]\nonumber\\&=&[[u\otimes t^m, v\otimes t^n],\a(w)\o t^k]+[[v\o
t^n,w\o t^k],\a(u)\o t^m]+[[w\o t^k,u\o t^m],\a(v)\o
t^n]\nonumber\\&=&\D_1\o t^{m+n+k}+\D_2\o t^{m+n+k-1}+\D_3\o
t^{m+n+k-2},
\end{eqnarray}
where
\begin{eqnarray}
\D_1&=&[[u,v],\a(w)]+[[v,w],\a(u)]+[[w,u],\a(v)],\label{de1}\\
\D_2
&=&m\big([u,v]\c\a(w)+[u\c
v,\a(w)]-\a(u)\c[v,w]-[u,w]\c\a(v)-[u\c w,\a(v)]\big)\nonumber\\
&&-n\big([v,u]\c\a(w)+[v\c u,\a(w)]-\a(v)\c[u,w]-[v,w]\c\a(u)-[v\c
w,\a(u)]\big)\nonumber\\&&+k\big(\a(w)\c[u,v]+[w,v]\c\a(u)+[w\c
v,\a(u)]-[w,u]\c\a(v)-[w\c u,\a(v)]\big),\label{de2}\\
\D_3&=&(m^2-m)\big((u\c v)\c\a(w)-(u\c
w)\c\a(v)\big)\nonumber\\&&+(n^2-n)\big((v\c w)\c\a(u)-(v\c
u)\c\a(w)\big)\nonumber\\
&&+(k^2-k)\big((w\c u)\c\a(v)-(w\c
v)\c\a(u)\big)\nonumber\\&&+mn\big((u\c v)\c\a(w)-(v\c
u)\c\a(w)-\a(u)\c(v\c w)+\a(v)\c(u\c w)\big)\nonumber\\&&+mk\big((w\c
u)\c\a(v)-(u\c w)\c\a(v)-\a(w)\c(u\c v)+\a(u)\c(w\c
v)\big)\nonumber\\&&+nk\big((v\c w)\c\a(u)-(w\c
v)\c\a(u)-\a(v)\c(w\c u)+\a(w)\c(v\c u)\big).\label{de3}
\end{eqnarray}
 By \eqref{JJ} and \eqref{de1}, we have $\D_1=0$. By \eqref{Hom-Nov1+}, \eqref{Hom-Nov2+} and \eqref{de3}, we get $\D_3=0$. By \eqref{HGB} and  \eqref{de2}, we have $\D_2=0$. Finally, \eqref{jj} gives
\begin{eqnarray}\label{jjj}
[[u\otimes t^m, v\otimes
t^n],\varphi(w\o t^k)]+[[v\o t^n,w\o t^k],\varphi(u\o
t^m)]+[[w\o t^k,u\o t^m],\varphi(v\o t^n)]=0.
\end{eqnarray}
Thus $(\L,[-,-],\varphi)$ is a Hom-Lie algebra.

The arguments above are reversible. Indeed, if $(\L,[-,-],\varphi)$
is a Hom-Lie algebra, then $\eqref{jjj}$ holds, namely, $\eqref{jj}=0$. This
forces $\D_1=0$ (so $(\mathcal
{A}, [\d,\d],\a)$ is a Hom-Lie algebra), and $\D_2=\D_3=0$ for all $m,n,k\in\mathbb{Z}$. In
particular, set $n=k=0$ and $m\neq 0, 1$ in \eqref{de2} and
\eqref{de3}.  In this case, $\D_2=0$ is equivalent to
$$[u,v]\c\a(w)+[u\c v,\a(w)]-\a(u)\c[v,w]-[u,w]\c\a(v)-[u\c
w,\a(v)]=0,$$ and $\D_3=0$ amounts to $$(u\c v)\c\a(w)=(u\c
w)\c\a(v).$$ Set $m=n=k=1$ in \eqref{de3}. Then $\D_3=0$
yields
$$(u\c v)\c\a(w)-(v\c u)\c\a(w)-\a(u)\c(v\c w)+\a(v)\c(u\c w)=0.$$
Because $u, v, w\in\A$ are arbitrary, we conclude that $(\mathcal
{A},[\d,\d],\circ,\a)$ is a Hom-GD bialgebra.  \QED\vskip8pt

Let $(\mathcal {A},[\d,\d],\circ,\a)$ be a Hom-GD bialgebra. For $u\in\A$, $m\in\mathbb{Z}$, set $u[m]=u\o t^m$. Then \eqref{aff} can be rewritten as
\begin{eqnarray}\label{re-aff}
[u[m], v[n]]=[u,v][m+n]+m (u\c v)[m+n-1]-n (v\c u)[m+n-1],
\end{eqnarray}
for all  $u, v\in\A$ and $m, n\in\mathbb{Z}$. Consider the following $\L$-valued formal distribution
\begin{eqnarray}\label{ge}
u(z)=\mbox{$\sum_{m\in\mathbb{Z}}$}u[m]z^{-m-1}\in\L[[z,z^{-1}]].
\end{eqnarray}
\eqref{re-aff} is equivalent to
\begin{eqnarray}\label{dis}
[u(z_1),v(z_2)]&=&[u, v](z_2)\delta(z_1-z_2)+\partial_{z_2}(u\c
v)(z_2)\delta(z_1- z_2)\nonumber\\&&+(u\c v+v\c
u)(z_2)\partial_{z_2}\delta(z_1-z_2).
\end{eqnarray}
Furthermore, we have
\begin{eqnarray}
(z_1-z_2)^2[u(z_1),v(z_2)]=0,
\end{eqnarray}
namely, $[u(z_1),v(z_2)]$ is local. Applying $F^\lambda_{z_1,z_2}$ to \eqref{dis}, we obtain
\begin{eqnarray}
F^\lambda_{z_1,z_2}([u(z_1),v(z_2)])=[u,v](z_2)+\partial_{z_2}(u\c v)(z_2)+\lambda(u\c v+v\c
u)(z_2),\label{fourier4}
\end{eqnarray}
from which we can define a $\l$-bracket
$[\d_{\l}\d]:\A\o\A\rightarrow \A[\lambda]$ by
\begin{eqnarray}\label{disss}
[u_{\lambda}v]=[u,v]+\partial(u\c v)+\lambda(u\c v+v\c u), \
\mbox{for}\ u,v\in\A,
\end{eqnarray}
where $\p=\p_{z_2}$. In terms $j$-products,
\eqref{disss} is equivalent to
\begin{eqnarray}\label{j-product2}
u_{(0)}v=[u,v]+\partial (u\c v),\ \ u_{(1)}v=u\c v+v\c u,\ \
u_{(j)}v=0, \ \mbox{for}\ u,v\in\A, \ j\geq 2,
\end{eqnarray}
which shows that $\A$ is closed under all the $j$-products. It is easy to check that the $\l$-bracket defined by \eqref{disss} satisfies \eqref{HL2}. However, the Hom-Jacobi identity \eqref{HL3} does not hold due to the Hom-Novikov algebra structure $(\A,\c)$. Indeed, as we will see it, the problem can be solved by changing the order of $u$ and $v$ on the right-hand side of \eqref{disss}.

 The following result is analogous to Gel'fand and Dorfman's statement for Gel'fand-Dorfman bialgebras and Xu's theorem for super Gel'fand-Dorfman
bialgebras.
\begin{theo}\label{main} Let $\mathcal {A}$ be a vector
space equipped with a linear endomorphism $\a$ and two operations $\circ$ and $[\d,\d]$. Let
$\mathcal{R}_{\A}=\mathbb{C}[\p]\o\A$ be the free $\C[\p]$-module over $\A$. Extend $\a$ to $\mathcal{R}_{\A}$ by
\begin{eqnarray}\label{*} \a(f(\p)\otimes u)=f(\p)\otimes \a(u), \ \mbox{for} \ u\in\A, \ f(\p)\in\C[\p].
\end{eqnarray}
and define a $\l$-bracket $[\d_\l
\d]:\A\otimes\A \rightarrow
\A[\lambda]$ by
\begin{eqnarray}\label{**}[u_{\lambda}v]=[v,u]+\partial(v\c
u)+\lambda(v\c u+u\c v), \ \mbox{for} \ u,v\in\A.
\end{eqnarray}
Then $(\mathcal{A},[\d,\d],\circ,\a)$ is a  Hom-GD bialgebra if and only if $(\mathcal{R}_\mathcal {A},[\d_\l\d],\a)$ is a
Hom-Lie conformal algebra of degree $2$.
\end{theo}
\noindent{\it Proof.~}Suppose that $(\mathcal {A},[\d,\d],\circ,\a)$
is a Hom-GD bialgebra. To show $(\mathcal{R}_\mathcal
{A},[\d_\l\d],\a)$ is a Hom-Lie conformal algebra, we need to
check the axioms from Definition \ref{HLCF}. The conformal
sesquilinearity \eqref{HL1} is naturally satisfied as the $\l$-bracket is defined on generators of a free $\C[\p]$-module, as it is explained in Remark \ref{remark}. It follows from \eqref{*} that $\a\p=\p\a$. Thus it
suffices to verify \eqref{HL2} and \eqref{HL3}.
 By \eqref{**} and \eqref{ss}, we have
\begin{eqnarray}\label{skew-s}
[v_{-\lambda-\p}u]&=&[u,v]+\partial(u\c v)+(-\lambda-\p)(v\c u+u\c
v)\nonumber\\&=&[u,v]-\partial(v\c u)-\lambda(v\c u+u\c v)\nonumber
\\&=&-[v,u]-\partial(v\c u)-\lambda(v\c u+u\c v)
=-[u_\l v],
\end{eqnarray}
which proves \eqref{HL2}. To show \eqref{HL3}, with a direct calculation we get
\begin{eqnarray*}
[\a(u)_\l [v_{\mu}
w]]&=&[[w,v],\a(u)]+\p[w,v]\c\a(u)+\l\big([w,v]\c\a(u)+\a(u)\c[w,v]\big)\\
&&+(\p+\l)\big([w\c v,\a(u)]+\p(w\c v)\c\a(u)+\l((w\c
v)\c\a(u)+\a(u)\c(w\c v))\big)\\
&&+\mu\big([v\c w+w\c v,\a(u)]+\p(v\c w+w\c v)\c\a(u)+\l(v\c w+w\c
v)\c\a(u)\big)\\&&+\mu\l\big(\a(u)\c(v\c w+w\c v)\big)\\
&=&[[w,v],\a(u)]+\p\big([w,v]\c\a(u)+[w\c v,\a(u)]\big)+\p^2(w\c v)\c\a(u)\\
&&+\l\big([w,v]\c\a(u)+\a(u)\c[w,v]+[w\c v,\a(u)]\big)+\mu\big([v\c w+w\c v,\a(u)]\big)\\
&&+\p\l\big(2(w\c v)\c\a(u)+\a(u)\c(w\c v)\big)+\l^2((w\c
v)\c\a(u)+\a(u)\c(w\c v)\big)\\
&&+\p\mu(v\c w+w\c v)\c\a(u)+\mu\l\big((v\c w+w\c
v)\c\a(u)+\a(u)\c(v\c w+w\c v)\big),
\end{eqnarray*}
where the second equality follows from \eqref{HL1} and \eqref{HL4}. Similarly, we
have
\\[4pt]\hspace*{13ex}$\dis [\a(v)_\mu [u_{\l} w]]
=[[w,u],\a(v)]+\p\big([w,u]\c\a(v)+[w\c u,\a(v)]\big)+\p^2(w\c u)\c\a(v)$
\\[4pt]\hspace*{13ex}$\dis\phantom{[\a(v)_\mu [u_{\l} w]]=}+\mu\big([w,u]\c\a(v)+\a(v)\c[w,u]+[w\c
u,\a(v)]\big)$
\\[4pt]\hspace*{13ex}$\dis\phantom{[\a(v)_\mu [u_{\l} w]]=}
+\l\big([u\c w+w\c u,\a(v)]\big)+\p\mu\big(2(w\c
u)\c\a(v)+\a(v)\c(w\c u)\big)$\\
\\[4pt]\hspace*{13ex}$\dis\phantom{[\a(v)_\mu [u_{\l} w]]=}
+\mu^2((w\c u)\c\a(v)+\a(v)\c(w\c
u))+\p\l(u\c w+w\c u)\c\a(v)$
\\[4pt]\hspace*{13ex}$\dis\phantom{[\a(v)_\mu [u_{\l} w]]=}
+\mu\l\big((u\c w+w\c
u)\c\a(v)+\a(v)\c(u\c w+w\c u)\big),$\\
\noindent and
\begin{eqnarray*}
[[u_{\l} v]_{\l+\mu}\a(w)]
&=&[\a(w),[v,u]]+\mu\big(\a(w)\c[v,u]+[v,u]\c\a(w)-[\a(w),v\c u]\big)\\
&&+\p\a(w)\c[v,u]+\l\big(\a(w)\c[v,u]+[v,u]\c\a(w)+[\a(w),u\c v]\big)\\
&&+\p\l\a(w)\c(u\c v)+\l^2\big(\a(w)\c(u\c v)+(u\c v)\c\a(w)\big)\\
&&-\p\mu\a(w)\c(v\c u)-\mu^2\big(\a(w)\c(v\c u)+(v\c u)\c\a(w)\big)\\
&&+\mu\l\big(\a(w)\c (u\c v)+(u\c v)\c\a(w)-\a(w)\c(v\c u)-(v\c
u)\c\a(w)\big).
\end{eqnarray*}
So \eqref{HL3} holds if and only if the following
equalities hold
\begin{eqnarray}
0&=&[[w,v],\a(u)]-[[w,u],\a(v)]-[\a(w),[v,u]],\label{ee1}\\
0&=&{[w,v]\c\a(u)+[w\c v,\a(u)]}-([w,u]\c\a(v)+[w\c
u,\a(v)]+\a(w)\c[v,u]),\label{ee2}\\
0&=&{(w\c v)\c\a(u)}-(w\c u)\c\a(v),\label{ee3}\\
0&=&{[w,v]\c\a(u)+\a(u)\c[w,v]+[w\c v,\a(u)]}-[u\c w+w\c
u,\a(v)]\nonumber\\&&-(\a(w)\c[v,u]+[v,u]\c\a(w)+[\a(w),u\c
v]),\label{ee4}\\
0&=&[v\c w+w\c v,\a(u)]-\big([w,u]\c\a(v)+\a(v)\c[w,u]+[w\c
u,\a(v)]\big)\nonumber\\&&-\big(\a(w)\c[v,u]+[v,u]\c\a(w)-[\a(w),v\c
u]\big),\label{ee5}\\0&=&\big(2(w\c v)\c\a(u)+\a(u)\c(w\c
v)\big)-(u\c w+w\c
u)\c\a(v)-\a(w)\c(u\c v),\label{ee6}\\
0&=&(v\c w+w\c v)\c\a(u)+\a(w)\c(v\c u)-\big(2(w\c
u)\c\a(v)+\a(v)\c(w\c u)\big),\label{ee7}\end{eqnarray}
\begin{eqnarray}
0&=&(w\c v)\c\a(u)+\a(u)\c(w\c v)-\a(w)\c(u\c v)-(u\c v)\c\a(w),\label{ee8}\\
0&=&(w\c u)\c\a(v)+\a(v)\c(w\c u)-\a(w)\c(v\c u)-(v\c u)\c\a(w),\label{ee9}\\
0&=&(v\c w+w\c v)\c\a(u)+\a(u)\c(v\c w+w\c v)-(u\c w+w\c
u)\c\a(v)\nonumber\\&&-\a(v)\c(u\c w+w\c u)-\a(w)\c (u\c v)-(u\c
v)\c\a(w)\nonumber\\&&+\a(w)\c(v\c u)+(v\c u)\c\a(w).\label{ee10}
\end{eqnarray}
Equalities \eqref{ee1}, \eqref{ee2} and \eqref{ee3} hold because of \eqref{JJ}, \eqref{HGB} and \eqref{Hom-Nov2+}, respectively. Since $u,v,w$ are arbitrary, equalities \eqref{ee4}, \eqref{ee6} and \eqref{ee8} are respectively equivalent to \eqref{ee5}, \eqref{ee7} and \eqref{ee9}. Furthermore, \eqref{ee4} can be rewritten as
\begin{eqnarray*}
&&\big([w\c v,\a(u)]+[w,v]\c\a(u)-\a(w)\c[v,u]-[w\c
u,\a(v)]-[w,u]\c\a(v)\big)\\&&\ \ \ +\big([u\c v],\a(w)]-[u\c
w,\a(v)]-[u,w]\c\a(v)+[u,v]\c\a(w)-\a(u)\c[v,w]=0,
\end{eqnarray*}
which is implied by \eqref{HGB}. By \eqref{Hom-Nov2+}, equality \eqref{ee6} amounts
to
\begin{eqnarray*}
(w\c u)\c\a(v)-\a(w)\c(u\c v)=(u\c w)\a(v)-\a(u)\c(w\c v),
\end{eqnarray*}
which holds due to \eqref{Hom-Nov1+}. Similarly, we get \eqref{ee8}.
To check \eqref{ee10}, we rewrite it as
\begin{eqnarray*}
&&\big((v\c w)\c\a(u)-\a(v)\c(w\c u)-(w\c v)\c\a(u)+\a(w)\c(v\c
u)\big)\\&&\ \ \ \ +\big((w\c u)\c\a(v)-\a(w)\c(u\c v)-(u\c
w)\c\a(v)+\a(u)\c(w\c v)\big)\\&&\ \ \ \ \ \ -\big((u\c
v)\c\a(w)-\a(u)\c(v\c w)-(v\c u)\c\a(w)+\a(v)\c(u\c w)\big)=0,
\end{eqnarray*}
which follows from \eqref{Hom-Nov1+}. Thus, \eqref{HL3} holds and $(\mathcal{R}_\mathcal {A},[\d_\l\d],\a)$ is a Hom-Lie
conformal algebra, which is of degree $2$ by \eqref{quadratic} and \eqref{**}.

Conversely, assume that $(\mathcal{R}_\mathcal {A},[\d_\l\d],\a)$ is a Hom-Lie conformal algebra. The Hom-Jacobi identity \eqref{HL3} gives \eqref{ee1}--\eqref{ee10}, from which we can deduce \eqref{Hom-Nov1+}--\eqref{JJ} and \eqref{HGB} by the discussions above. Therefore, $(\mathcal {A},[\d,\d],\circ,\a)$ is a Hom-GD bialgebra.
 \QED

 \vskip10pt

\small\noindent{\bf Acknowledgements~}~{\footnotesize This work was
conducted during the author's visit to Rutgers University. She would
like to thank Professor Haisheng Li and his family for the hospitality during her stay there. Special thinks to Professor
Haisheng Li for comments on this work and valuable discussions on
vertex algebras. This work was supported by National Natural Science
Foundation grants of China (11301109), the Research Fund for the Doctoral Program of Higher Education
(20132302120042), China Postdoctoral Science
Foundation funded project (2012M520715), the Fundamental Research
Funds for the Central Universities (HIT. NSRIF. 201462), a grant from China Scholarship Council
(201206125047) and partially by National Natural Science
Foundation grants of China (11101387).}\\

\end{document}